\documentclass[12pt]{article}
\usepackage{amsmath}
\usepackage{amsthm}
\usepackage{amsfonts}
\newcommand{\Rset}{\mathbb{R}}
\newcommand{\Cset}{\mathbb{C}}

\newtheorem{thm}{THEOREM}

\newtheorem{prop}{PROPOSITION}
\newtheorem{cor}{COROLLARY}
\theoremstyle{remark}
\newtheorem{rem}{REMARK}
\begin{document}

\title{A New Factorization Property of the Selfdecomposable Probability
Measures}
\date{April 18, 2003}
\author{Aleksander M. Iksanov (Kiev), Zbigniew J. Jurek\footnote{This work
was completed while this author was visiting Wayne State University.}
(Wroc\l aw)\\
and Bertram M. Schreiber (Detroit)} \maketitle

\begin{quote}
\noindent {\footnotesize \textbf{ABSTRACT.} We prove that the convolution
of a selfdecomposable distribution with its background driving law is again
selfdecomposable if and only if the background driving law is
s-selfdecomposable. We will refer to this as the \textit{factorization
property} of a selfdecomposable distribution; let $L^f$ denote the set of
all these distributions. The algebraic structure and various
characterizations of $L^f$ are studied. Some examples are discussed, the
most interesting one
 being given by the L\'evy stochastic area integral. A nested
family of subclasses $L^{f}_n , n\ge 0,$ (or a filtration) of the class
$L^f$ is given.
\medskip

MSC\,2000 \textit{subject classifications.} Primary 60E07, 60B12; secondary
60G51, 60H05.

}
\end{quote}

\medskip
Abbreviated title:
\begin{center}
 SELFDECOMPOSABLE MEASURES
\end{center}

\vfill\pagebreak Limit distribution theory and the study of infinitely
divisible distributions belong to the core of probability and mathematical
statistics. Here we investigate an unexpected relation between two classes
of distributions, $L$ and $\mathcal{U}$, each of which can be defined in
terms of a collection of inequalities involving the L\'evy spectral measure
and a semigroup of mappings. Indeed, the L\'evy class $L$ is defined via
\textit{linear transformations} while the class $\mathcal{U}$ involves
\textit{nonlinear transformations}. Yet these classes exhibit some
similarities and relationships, such as the proper inclusion $L \subset
\mathcal{U}$; cf. Jurek (1985).

In recent years class $L$ distributions have found many applications, in
particular through their BDLP's (background driving L\'evy processes); cf.
for example Barndorff-Nielsen and Shephard (2001) and the references there.
Also there were developed stochastic methods for finding the BDLP's of some
selfdecomposable distributions; cf. Jeanblanc, Pitman and Yor (2002). On
the other hand, in Jacod, Jakubowski and M\'emin (2003), class
$\mathcal{U}$ distributions appeared in the context of an approximation of
processes by their discretization.

In Section 1 we recall the definitions of the class $L$ of
selfdecomposable distributions and the class $\mathcal{U}$ of
s-selfdecomposable ones, followed by their random integral
representations. In Section 2 we introduce the new notion termed
the \textit{factorization property} and the corresponding class
$L^f$. These are class $L$ distributions whose convolutions with
their background driving distributions are again class $L$
distribution. Elements of the class $L^f$ are characterized in
terms of their Fourier transforms (Corollary 1 and Theorem 3) and
their L\'evy spectral measures (Corollary 2 and Theorem 2).
Proposition 1 describes the topological and algebraic structures
of the class $L^f$. Some explicit examples of $L^f$ probability
distributions, which includes among others the L\'evy stochastic
area integral, are given in Section 3.

Our main results are given in the generality of probability
measures on a Banach space, but they are new for distributions on
the real line as well. Indeed, changing the pairing between a
Banach space and its dual to the scalar product in all our proofs,
one gets results in Euclidean spaces and Hilbert spaces. On the
other hand, if one deals with variables assuming values in
function spaces (stochastic processes), then Banach spaces provide
the natural setting. For instance, Brownian motion and Bessel
processes, when restricted to finite time, can be viewed as
measures on Banach spaces of continuous functions. There is a vast
literature dealing with probability on Banach spaces (e.g., cf.
Araujo-Gine (1980) or \cite{P} and the references in the articles
found there), and much of the work leading up to the results
presented here was carried out in this context. Finally, our
proofs do not depend on the dimension of the space on which the
probability measures are defined. Thus the generality of Banach
spaces seems to be the natural one. This paper continues the
investigations of Jurek (1985).

\medskip
\textbf{1. Introduction and notation.} Let $E$ denotes a real separable
Banach space, $E^{\prime}$ its conjugate space, $< \cdot,\cdot >$ the usual
pairing between $E$ and $E^{\prime}$, and $||.||$ the norm on E. The
$\sigma$-field of all Borel subsets of $E$ is denoted by $\mathcal{B}$,
while $\mathcal{B}_{0}$ denotes Borel subsets of $E \setminus \{0\}$. By
$\mathcal{P}(E)$ we denote the (topological) semigroup of all Borel
probability measures on $E$, with convolution ``$\ast$" and the topology of
weak convergence ``$\Rightarrow$." As in Jurek (1985), we denote the closed
subsemigroup of infinitely divisible measures in $\mathcal{P}(E)$ by
$ID(E)$.

Each ID distribution $\mu$ is uniquely determined by a triple: a shift
vector $a\in E$, a Gaussian covariance operator $R$, and a L\'evy spectral
measure $M$; we will write $\mu = [a,R,M]$. These are the parameters in the
L\'evy-Khintchine representation of the characteristic function
$\hat{\mu}$, namely $\mu \in ID$ iff $\hat{\mu}(y)= \exp(\Phi(y))$, where
\begin{align*}
\Phi(y) &= i<y,a>-1/2<Ry,y>+ \nonumber \\
  &\quad \int_{E \setminus\{0\}} [e^{i<y,x>}-1-i<y,x>1_{||x||\leq 1}(x)]M(dx),\ \ y
\in E^{\prime};
\end{align*}
$\Phi$ is called the \textit{L\'evy exponent} of $\hat{\mu}$ (cf. Araujo
and Gin\'e (1980), Section 3.6).

On the Banach space $E$ we define two families of transforms $T_{r}$ and
$U_{r}$, for $r>0$, as follows:
\[
T_{r}x = rx  \ \ \mbox{and} \ \ U_{r}x = \max(0,||x||-r)\frac{x} {||x||}, \
\ U_{r}(0) = 0.
\]
The $T_r$'s are linear mappings; the $U_r$'s are nonlinear and are called
\textit{shrinking operations}, or \textit{s-operations} for short.

In Jurek (1985) the class $L(E)$ of \textit{selfdecomposable} measures was
introduced as those $\mu = [a,R,M]\in ID(E)$ such that
\[
M \ge T_cM, \text{  for  } 0 < c < 1.
\]
As pointed out there (Corollary 3.3), this condition is equivalent to the
traditional definition:
\begin{equation}
\mu \in L(E) \ \ \mbox{iff} \ \ \forall\,{(0<c<1)}\ \exists\,(\mu_{c}\in
\mathcal{P}(E)) \ \ \mu =T_{c}\mu \ast \mu_{c}.
\end{equation}
It follows easily that $L(E)$ is a closed convolution topological
semigroup of $\mathcal{P}(E)$. The importance of the class L(E)
arises from the fact that it extends the classical and
much-studied class of \textit{the stable distributions}.

One important example of a selfdecomposable measure is Wiener
measure $\mathcal{W}$ on the Banach space $C_{\mathbb{R}}([0,1])$.
That $\mathcal{W}$ is selfdecomposable follows immediately from
the fact that its finite-dimensional projections are Gaussian
measures, hence selfdecomposable.

Similarly, a measure $\mu = [a,R,M]$ is called \textit{s-selfdecomposable}
on $E$, and we will write $\nu \in \mathcal{U}(E)$, if
\[
M \ge U_rM, \text{  for  } 0 < r < \infty.
\]
As shown in Jurek (1985) (cf. Jurek (1981)),
\begin{equation}
\nu \in \mathcal{U}(E)\ \ \mbox{iff}\ \ \forall\,{(0<c<1)}\ \exists \,
(\nu_{c}\in\mathcal{P}(E))  \ \ \nu =T_{c}\nu ^{\ast c}\ast \nu_{c}
\end{equation}
(the convolution power is well defined as $\nu$ is in $ID(E)$). In
particular, we infer that $(\mathcal{U}(E)$ is also a closed convolution
topological semigroup. In fact, we have the inclusions
\[
L(E)\subset \mathcal{U}(E)\subset ID(E) \subset \mathcal{P}(E).
\]

\medskip
Relations between the semigroups $L(E)$ and $\mathcal{U}(E)$ and their
characterizations were studied in Jurek (1985), and this paper is our main
reference for this work, including the terminology and basic notation.

Let
\[
ID_{log}(E) = \{\mu\in ID(E): \int_E \log(1 + \|x\|)\,\mu(dx) <
\infty \},
\]
and recall that the mapping $\mathcal{I}:ID_{log}(E) \to L(E)$ given by
\begin{equation}
\mathcal{I}(\rho) = \mathcal{L}(\int_{(0,\infty)}e^{-s}dY_\rho(s)),
\end{equation}
is an algebraic isomorphism between the convolution semigroups
$ID_{log}$ and $L$; cf. Jurek (1985), Theorem 3.6. Above
$Y_\rho(\cdot)$ denotes a L\'evy process, i.e., an $E$-valued
process with stationary and independent increments, with
trajectories in the Skorohod space of cadlag functions, and such
that $Y_\rho(0)=0$ a.s and $\mathcal{L}(Y_\rho(1))=\rho$.

If a class $L$ distribution is given by (3) then we  refer to $Y_{\rho}$ as
the \textit{\textbf{b}ackground \textbf{d}riving \textbf{L}\'evy
\textbf{p}rocess} (in short, the BDLP, cf. Jurek (1996)). The measure
$\rho$ in (3) will be called the \textit{\textbf{b}ackground
\textbf{d}riving \textbf{p}robability \textbf{d}istribution} (in short, the
BDPD) and the r.v. $Y_{\rho}(1)$ is the \textit{\textbf{b}ackground
\textbf{d}riving \textbf{r}andom \textbf{v}ariable} (in short, the BDRV).

Similarly, for s-selfdecomposable distributions we define a mapping
$\mathcal{J}: ID(E) \to \mathcal{U}(E)$ given by
\begin{equation}
\mathcal{J}(\rho) = \mathcal{L}(\int_{(0,1)}s\, dY_\rho(s)),
\end{equation}
which is an isomorphism between the topological semigroups $ID(E)$
and $\mathcal{U}(E)$; cf. Jurek (1985), Theorem 2.6. In (4)
$Y_{\rho}(\cdot)$ is an arbitrary  L\'evy process.

\medskip
Let $\hat{\mu}(y) = \int_E e^{i<y,x>}\mu(dx),\ \ y \in E^{\prime}$, be
\textit{the characteristic function} (Fourier transform) of a measure
$\mu$. Then random integrals like (3) or (4) have characteristic functions
of the form
\begin{equation}
\left(\mathcal{L}(\int_{(a,b]}h(t)dY_{\rho}(t))\right)^{\widehat{ }}(y)=
\exp \int_{(a,b]}\log \hat{\rho}(h(t)y)dt,
\end{equation}
when $h$ is a deterministic function and $Y_{\rho}(.)$ a L\'evy process;
cf. Lemma 1.1 in Jurek (1985).

\medskip
\textbf{2. A new factorization property of class $L$
distributions.} Class $L$ distributions decompose by themselves as
is evident from the convolution equation (1). However, recently it
has been noted that in some  classical formulae class L
distributions appear convoluted with their background driving
probability distributions (BDPD); for instance the L\'evy
stochastic area integral is one such example; cf. Jurek (2001).
The following is our main result that describes the cases when a
selfdecomposable distribution can be factored as another class L
distribution and its corresponding BDPD.
\medskip
\begin{thm}
A selfdecomposable probability distribution $\mu=\mathcal{I}(\nu)$
convoluted with its background driving law $\nu$ is
selfdecomposable if and only if $\nu$ is s-selfdecomposable. More
explicitly, for $\nu$ and $\rho$ in $ID_{log}$ we have
\begin{equation}
\mathcal{I}(\nu)\ast \nu = \mathcal{I}(\rho) \ \ \mbox{iff} \ \
\nu=\mathcal{J}(\rho).
\end{equation}
\end{thm}

\textit{Proof. Sufficiency.} Suppose $\mu$ is selfdecomposable and its
background driving law $\nu$ is s-selfdecomposable. That is
$\mu=\mathcal{I}(\nu)$ for some unique $\nu \in ID_{log}$ and
$\nu=\mathcal{J}(\rho)$. Hence $\rho \in ID_{log}(E)$, by formula (4.1) in
Jurek (1985). Consequently we have
\[
\nu \ast \mathcal{I}(\nu)= \mathcal{J}(\rho) \ast
\mathcal{I}(\mathcal{J}(\rho))= \mathcal{J}(\rho \ast
\mathcal{I}(\rho))=\mathcal{I}(\rho)\in L,
\]
where the last equality follows from Corollary 4.6 in Jurek (1985). Also we
have used  the fact that the mappings $\mathcal{I}$ and $\mathcal{J}$
commute; cf. \textit{ibid.}, Theorem 3.6 and Corollary 4.2. The sufficiency
is proved.

\textit{Necessity.} Suppose that a selfdecomposable $\mu=\mathcal{I}(\nu)$
is such that $\nu \ast \mathcal{I}(\nu)$ is again selfdecomposable. Then
there is a unique $\rho \in ID_{log}$ such that
\[
\nu \ast \mathcal{I}(\nu)= \mathcal{I}(\rho).
\]
Applying the mapping $\mathcal{J}$ to both sides and employing Corollary
4.6 in Jurek (1985) and the commutativity, we conclude
\[
\mathcal{I}(\nu)=\mathcal{J}(\nu\ast\mathcal{I}(\nu))=
\mathcal{J}(\mathcal{I}(\rho))=\mathcal{I}(\mathcal{J}(\rho)).
\]
Since $\mathcal{I}$ is one-to-one, $\nu=\mathcal{J}(\rho)$, which completes
the proof of necessity.

\medskip
We will say that a selfdecomposable probability distribution $\mu$ has the
\textit{factorization property} (we will write $\mu \in L^f$), if its
convolution with its BDPD gives another selfdecomposable distribution,
i.e.,
\begin{equation}
\mu \in L^f \ \ \ \ \mbox{iff} \ \ \mu=\mathcal{I}(\nu) \ \
\mbox{\emph{for} $\nu \in ID_{log}$ \emph{and}} \ \ \mu \ast \nu
\in L.
\end{equation}

\medskip
Before describing the algebraic structure of the class $L^f$, let us recall
that by definition \textit{class $L_1$ distributions} are those
selfdecomposable distributions for which the cofactors $\mu_c$ in (1) are
in $L(E)$. Equivalently, these are distributions of random integrals (3),
where $\rho$ is a distribution from the class $L$. This class was first
introduced for real-valued r.v.'s in Urbanik(1973) as the first of a
decreasing sequence $L_n$ ($ n = 0,1,2,...$) of subclasses of the class $L$
and later studied by Kumar and Schreiber (1978) and in the vector-valued
case in Kumar and Schreiber (1979), Sato (1980), and Jurek (1983a,b). In
fact, Jurek (1983a) contains the most general setting, where in (1) the
operators $T_a$ may be chosen from any one-parameter group of operators.

\begin{prop}
The class $L^f$ of selfdecomposable distributions with the
factorization property is a closed convolution subsemigroup of
$L$. Moreover:
\begin{itemize}
\item[\em{(i)}] For $a>0$, $T_a \mu \in L^f$ iff $\mu \in L^f$.

\item[\em{(ii)}] A probability measure $\mu \in L^f$ iff there
exists a (unique) probability measure $\nu \in ID_{log}$ such that
\begin{equation}
\mu=\mathcal{I}(\mathcal{J}(\nu))=\mathcal{J}(\mathcal{I}(\nu)), \
\ \mbox{i.e.,} \ \
L^f=\mathcal{I}(\mathcal{J}(ID_{log}))=\mathcal{J}(L).
\end{equation}

\item[\em{(iii)}] $L_1 \subset L^f$, where $L_1$ consists of those class
$L$ distributions whose $BDLP$ are in class $L$.
\end{itemize}
\end{prop}

\textit{Proof.} The semigroup structure and properties (i) and (ii) follow
from formula (6), Theorem 1, and properties of the mappings $\mathcal{I}$
and $\mathcal{J}$. To prove that $L^f$ is closed, let $\mu_n \in L^f$ and
let $\mu_n \Rightarrow \mu$. Then $\mu =\mathcal{I}(\nu) \in L$ by (3), and
$\mu_n=\mathcal{I}(\nu_n)\Rightarrow \mathcal{I}(\nu)=\mu$. From Jurek and
Rosinski (1988) we conclude that $\nu_n \Rightarrow \nu$ and
$\int_{E}\log(1+||x||)\nu_{n}(dx)\to \int_{E}\log(1+||x||)\nu(dx)$.
Consequently $\mu \ast \nu \in L$ and therefore $\mu\in L^f$, which proves
that $L^f$ is closed.

Since each $\mu \in L_1$ has its BDPD $\nu \in L$ and $L_1 \subset{L}$,
(iii) follows from the semigroup property of $L$.

\medskip
Theorem 1 can be expressed in terms of characteristic functions. Namely:

\begin{cor}
In order that
\begin{equation}
\exp\left(\int_{0}^{\infty}\log\hat{\nu}(e^{-s}y)ds\right) \cdot
\hat{\nu}(y)= \exp\int_{0}^{\infty}\log\hat{\rho}(e^{-s}y)ds,\ \ y \in E',
\end{equation}
for some $\nu$ and $\rho$ in $ID_{log}$, it is necessary and
sufficient that
\begin{equation}
\hat{\nu}(y)= \exp \int_{0}^{1} \log \hat{\rho}(sy)ds
\end{equation}
\end{cor}
The above follows from (5) and (6). For details cf. Jurek (1985), Theorems
2.9 and 3.10.

\begin{cor}
In order to have the equality
\[
\int_{(0,\infty)}N(e^sA)ds+N(A) = \int_{(0,\infty)}G(e^sA)ds, \ \
\mbox{for all}\ \ A \in \mathcal{B}_0 ,
\]
for some L\'evy spectral measures N and G with finite logarithmic moments
on sets $\{x:||x|| > c\}$, it is necessary and sufficient that
\begin{equation}
N(A)=\int_{(0,1)}G(t^{-1}A)dt, \ \ \mbox{for all} \ \ A \in
\mathcal{B}_0 .
\end{equation}
\end{cor}
This is easily obtained from (7), (9) and (10). For more details
cf. Jurek (1985) formulae (2.9) and (3.4).

\medskip
One may also characterize the factorization property purely in
terms of L\'evy spectral measures and functions, as shift and
Gaussian parts do not contribute any restrictions. For that
purpose let us recall that by the \textit{L\'evy spectral
function} of $\mu=[a,R,M]$ we mean the function
\begin{equation}
L_{M}(D,r): = - M(\{x \in E : \|x\|>r \ \text{and} \ x/||x|| \in D \}),
\end{equation}
where $D$ is a Borel subset of unit sphere $S = \{x\in E:
||x||=1\}$ and $r>0$. Note that $L_M$ uniquely determines $M$.

\begin{thm}
In order that $\mu=[a,R,M]$ have the factorization property, i.e., $\mu \in
L^f$, it is necessary and sufficient that there exist a unique L\'evy
spectral measure $G$ with finite logarithmic moments on all sets of form
$\{x:||x||>c\}$, $c>0$, such that
\begin{equation}
M(A)=\int_0^\infty \int_0^1 G(e^{t} s^{-1}A)ds\,dt, \ \ \mbox{for all} \ \
A \in \mathcal{B}_0 .
\end{equation}
Equivalently, for all Borel subsets $D$ of the unit sphere in $E$,
$dL_M(D,\cdot)/dr$ exists and
\begin{equation}
r \mapsto r \frac{dL_M(D,r)}{dr}
\end{equation}
is a convex, nonincreasing function on $(0,\infty)$.
\end{thm}

\textit{Proof.} If $\mu\in L^f$, then since $M(A) =
\int_0^{\infty} N(e^s A)ds$ and $N$ has the form (11), we infer
equality (13). From (13) we get
\begin{equation}
L_{M}(D,r)=\int_r^\infty \int_u^\infty \frac{L_{G}(D,w)}{w^2}dw du
= \int_r^\infty\frac{w-r}{w^2}L_{G}(D,w)dw,
\end{equation}
and consequently
\begin{equation}
\frac{d}{dr}\Big(r\frac{dL_{M}(D,r)}{dr}\Big)= - \int_r^\infty
\frac{dL_G(D,w)}{w},
\end{equation}
at points of continuity of $L_G(D,\cdot)$. Hence the existence of the first
derivative and the properties of the function (14) follow.

Conversely, if the function (14) is nonincreasing and convex, then
first of all, the L\'evy spectral measure $M$ corresponds to a
class $L$ probability measure, say $\mu$, by Jurek (1985), Theorem
3.2 (b). Furthermore, $\mu$ is of the form (3), where the BDRV
$Y_{\rho}(1)$ has finite logarithmic moment, and its L\'evy
spectral measure $G$ satisfies
\[
M(A) = \int_0^{\infty} G(e^s A)ds,\qquad A\in \mathcal{B}_0.
\]
Hence, in terms of the corresponding L\'evy spectral functions, the
convexity assumption implies that
\[
L_G(D,r) = -r \frac{dL_M(D,r)}{dr} = -\int_r^{\infty} q(D,s)ds,
\]
for a uniquely determined, nonincreasing, right-continuous
function $q(D,\cdot)$. In other words, $dL_G(D,r)/dr = q(D,r)$
exists almost everywhere and is nonincreasing in $r$. By Theorem
2.2 (b) in Jurek (1985) we infer that $G$ corresponds to a class
$\mathcal{U}$ probability measure, meaning that the distribution
of $Y_{\rho}(1)$ in (3) is in $\mathcal{U}$. By Theorem 1 we
conclude that $\mu$ has the factorization property, which
completes the proof.

\medskip
As an immediate consequence of (14) we have the following.

\begin{cor}
If a L\'evy spectral function $L_M(D,r)$ is twice differentiable
in $r$, then $\mu = [a,R,M]\in L^f$ if and only if the functions
$r\mapsto r^2 d^2L_M(D,r)/dr^2$ are nondecreasind and right
continuous on $(0,\infty)$.
\end{cor}

\medskip
Finally we describe distributions in the class $L^f$ in terms of their
characteristic functions.

\begin{thm}
\emph{(a)} A function $g :E^{\prime}\to \Cset$ is the
characteristic function of a class $L^f$ distribution if and only
if there exists a unique $\nu \in ID_{log}$ such that
\begin{equation}
g(y) = \exp \left[\int_0^1\int_0^w
\frac{\log\hat{\nu}(uy)}{w^2}du\,dw \right] =
\frac{\mathcal{I}(\nu)^{\widehat{\ }}
(y)}{\mathcal{J}(\nu)^{\widehat{\ }}(y)}, \ \ y \in E^{\prime}.
\end{equation}

\emph{(b)} $\Phi(\cdot)$ is the L\'evy exponent of a class $L^f$
distribution if and only if for each $y \in E^{\prime}$, the function
$\Rset \ni t \mapsto \Phi(ty)\in \Cset$ is twice differentiable and
\[
\Psi(y) = \left.\left[2\frac{d}{dt}\Phi(ty) +
\frac{d^2}{dt^2}\Phi(ty)\right]\right|_{t=1}
\]
is the L\'evy exponent of a distribution from $ID_{log}$.
\end{thm}

\textit{Proof.} (a) If $\mu \in L^f$ then by (8), (3)--(5), changing
variables and the order of integration gives
\begin{align*}
\log\hat{\mu}(y) &= \int_0^\infty
\log(\mathcal{J}(\nu)^{\widehat{\ }})(e^{-s}y)ds =
\int_0^{\infty}\int_0^1\log \hat{\nu}(e^{-s}ty)dt\,ds \\
  &= \int_0^\infty \int_0^{e^{-s}}\log\hat{\nu}(uy)e^sdu\,ds
  = \int_0^1\int_0^w \frac{\log\hat{\nu}(uy)}{w^2}du\,dw.
\end{align*}
The other equality follows from Theorem 1.

Conversely, if $\nu \in ID_{log}$ then the random integral
$\mathcal{I}(\nu)$ exists and consequently $\mu =
\mathcal{J}(\mathcal{I}(\nu))$ is defined as well. So the
calculation above applies to show that its characteristic function
is of the form (17), which completes proof of part (a).

For part (b) note that if $\mu$ and $\nu$ are related as above, with
respective L\'evy exponents $\Phi$ and $\Psi$, then by (a)
\begin{align}
r_y(t) &= \Phi(ty) = \int_0^1\int_0^w \frac{\Psi(tuy)}{w^2}du\,dw \nonumber \\
   & = \int_0^1 \int_0^{tw} \frac{\Psi(vy)}{tw^2}dv\,dw = \int_0^t\int_0^s
   \frac{\Psi(vy)}{s^2}dv\,ds.
\end{align}
Hence $r_y$ is twice differentiable, and the formula for $\Psi$ follows.

On the other hand, since $r_{sy}(t) = r_y(st)$, we have $r_{sy}'(t) =
sr_y'(st)$ and $r_{sy}''(t) = s^2r_y''(t)$. If $\Phi$ and $\Psi$ are
related as in (b), then
\[
\Psi(sy) = 2r_{sy}'(1) + r_{sy}''(1) = 2sr_y'(s) + s^2r_y''(s) =
\frac{d}{ds} \left[s^2 \frac{d}{ds}\Phi(sy)\right].
\]
Since $\Psi(0) = 0$, two integrations give (18), so by (a) the proof is
complete.

\medskip
As was already  mentioned above, K. Urbanik (1973) introduced a family of
decreasing classes of \textit{n times selfdecomposable distributions}
\begin{equation}
ID \supset L_0 \supset L_1 \supset ... \supset L_n \supset L_{n+1}...
\supset L_{\infty} = \bigcap_{n=1}^{\infty}L_n \supset S ,
\end{equation}
via some limiting procedures, where $L_{0} \equiv L$ is the class of all
selfdecomposable distributions and $S$ denotes the class of all stable
distributions (in the above inclusions we suppressed the dependence of
classes $L_n$ on the Banach space $E$). Also note that the class $L_1$ of
distributions in Proposition 1 is exactly the class $L_1$ in the sequence
(19).

Recall that $\mu$ is $n$ times selfdecomposable if and only if it admits
the integral representation (3) with $\rho$ being $(n-1)$ times
selfdecomposable. For other equivalent approaches see Kumar and Schreiber
(1979), Sato (1980), or Jurek (1983a,b). Let us define classes $L^{f}_n$ of
measures with the \textit{class $L_n$ factorization property} as follows:
\begin{equation}
L^{f}_n = \{\mu \in L_n : \mu\ast \mathcal{I}^{-1}(\mu) \in L_n \}, \quad n
= 0,1,2,\ldots ,
\end{equation}
where the isomorphism $\mathcal{I}$ is given by (3) and
$\mathcal{I}^{-1}(\mu)$ is the BDPD (the probability distribution of the
BDRV Y(1)) for $\mu$. In other words, $\mu$ from $L_n$ is in $L^{f}_n$ if
when it is convolved with its BDPD one obtains another distribution from
the class $L_n$.

For the purpose of the next results let us recall that
\[
\mu \in L_n \ \ \mbox{iff} \ \ \mu= \mathcal{I}(\rho) \ \
\mbox{for a unique} \ \  \rho \in L_{n-1} \cap ID_{log};\ \ \
L_n=\mathcal{I}(L_{n-1} \cap ID_{\log}),
\]
and
\begin{multline*}
\mu \in L_n \ \ \mbox{iff} \\
\mu = \mathcal{L}\left(\int_{0}^{\infty}
e^{-s}dY_{\nu}\left(\frac{s^{n+1}}{(n+1)!}\right)\right) =
\mathcal{I}^{n+1}(\nu) \ \ \mbox{for a unique} \ \ \nu \in ID_{log^{n+1}}
\end{multline*}
(cf. Jurek (1983b).

\begin{prop}
For $n = 0,1,2,...$, we have that:
\begin{itemize}
\item[\emph{(i)}] The classes $L^{f}_n$ are closed convolution
semigroups also closed under the dilations $T_a,\ a>0$.

\item[\emph{(ii)}] $L_{n+1} \subset L^{f}_n \subset L_n$ \ \ \mbox{(proper
inclusions).}

\item[\emph{(iii)}] A probability measure $\mu \in L^{f}_n$ iff
there exists a unique probability measure $\nu \in ID_{log^{n+1}} = \{\rho
\in ID : \int_E \log^{n+1}(1+||x||)\rho(dx)< \infty \}$ such that $\mu =
\mathcal{I}^{n+1}(\mathcal{J}(\nu))$, where $\mathcal{I}^{1} = \mathcal{I}$
and $\mathcal{I}^{n}(\cdot) = \mathcal{I}(\mathcal{I}^{n-1}(\cdot))$ (i.e.,
the mapping $\mathcal{I}$ is composed with itself $n$ times). That is,
\begin{equation}
L^{f}_n = \mathcal{J}(L_{n})= \mathcal{I}(L^{f}_{n-1}\cap ID_{\log}),\ n\ge
0,\ \mbox{where} \ \ L^f_{-1} = \mathcal{J}(ID_{\log}).
\end{equation}
\end{itemize}
\end{prop}

\textit{Proof.} Part (i) follows the proof of Proposition 1 and definition
(20). Part (ii), for $n = 0$, is just Proposition 1 (ii). Suppose the
inclusions in (ii) hold for some $k \ge 1$. Then
\[
\mathcal{I}(L_{k+1} \cap ID_{\log}) \subset \mathcal{I}(L^{f}_k
\cap ID_{\log}) \subset \mathcal{I}(L_k \cap ID_{\log}),
\]
which means that  $L_{k+2} \subset L^{f}_{k+1} \subset L_{k+1}$,
and  therefore (ii) is proved for all $k$.

Since the mappings $\mathcal{I}$ and  $\mathcal{J}$ are one-to-one, to
prove that the inclusions are proper it suffices to notice that $L$ is a
proper subset of $\mathcal{U}$. The latter is true because in order for an
s-selfdecoposable $\mathcal{J}(\rho),\,\rho\in ID$, to be selfdecomposable,
i.e, equal to $\mathcal{I}(\nu)$ for some $\nu \in ID_{\log}$, it is
necessary and sufficient that $\rho= \mathcal{I}(\nu)\ast \nu$; cf. Jurek
(1985), Theorem 4.5.

For part (iii) we again use induction argument, Proposition 1, and
the characterization of the classes $L_n$ quoted before
Proposition 2.

\medskip
\begin{rem}
From formula (21) in Proposition 2, we see that the sequence of
classes $L^{f}_n$ is obtained from the sequence $L_n$ in (19) by
applying the mapping $\mathcal{J}$ and then inserting it to the
right. This produces the following sequence of interlacing
subclasses:
\begin{multline}
\mathcal{U} \supset L_0 \supset L^{f}_0 \supset L_1 \supset L^{f}_1 \supset L_2 \supset ... \\
\supset L_n \supset L^{f}_n \supset L_{n+1}\supset ... \supset L_{\infty}=
\bigcap_{n=1}^{\infty}L_n = \bigcap_{n=1}^{\infty}L^{f}_n \supset S .
\end{multline}
For the last inclusion recall that stable distributions are in $L$ and have
stable laws as their BDPD. In other words, the stable distributions are
invariant under the mapping $\mathcal{I}$. In fact, the same is true for
$\mathcal{J}$ (cf. Jurek (1985), Theorems 3.9 and 2.8).
\end{rem}

\medskip
\begin{rem}
Using Proposition 2 (iii) inductively one can obtain characterizations of
the classes $L_n^f,\ n\ge 1,$ similar to those obtained for $L_0^f \equiv
L^f$ in Corollaries 3 and 4 or in Theorem 3. We leave these calculations
for the interested reader.
\end{rem}

\medskip
\textbf{3. Examples of distributions with the factorization
property.} In this last section we provide some explicit examples
both on arbitrary Banach space and on the real line.

\textbf{A.} 1. On any Banach space, as pointed out in Remark 1,
all stable measures have the factorization property. In fact, a
stable measure has a stable processes as its BDLP; cf. also Jurek
(1985), Theorem 3.8.

2. On any Banach space $E$, for positive constants $\alpha, \beta$
and  vector $z$ on the unit sphere in $E$ let
\[
K_{\alpha, \beta, z}(A)= \alpha
\int_0^{\beta}(\frac{\beta}{v}-1)\delta_{vz}(A)dv, \quad A \in
\mathcal{B}_0.
\]
Then the infinitely divisible measures $[a, 0, K_{\alpha, \beta,
z}]$, $a \in E$, have the factorization property. This follows by
applying Theorem 2 with $G=\,c\delta_x$, for $c>0$ and $0\neq x\in
E$; explicitly $\alpha = c/||x||, \, \beta=||x||,\, z=x/||x||$.

Because of Proposition 1, dilations, convolutions and weak limits
of the above probability measures possess the factorization
property as well; also cf. Jurek (1985), Theorem 2.10.

\textbf{B.} For our examples on real line we need some auxiliary
facts.

Recall that for an $ID$ distribution on the real line with L\'evy
spectral measure $M$, its L\'evy spectral function, as defined
above, is separately given on the positive and negative half-lines
as follows:
\begin{equation}
L_{M}(x) =
\begin{cases}
- M(\{s\in \Rset | s >x\}, \quad &\mbox{for}\ x>0 \\
  M\{s\in \Rset| s<x \}, &\mbox{for}\ x<0
\end{cases} .
\end{equation}
Indeed, for $x > 0$ we set $L_M(x) = L_M(\{1\},x)$, while for $x < 0$, set
$L_M(x) = -L_M(\{-1\},x)$.

From O'Connor (1979) or Jurek (1985),Theorem 2.2 (b) we have the following
description of s-selfdecomposable distributions :
\begin{align}
[a,\sigma^2,M] \in \mathcal{U}(\Rset)\ \mbox{iff}\ &L_{M}(x)\ \mbox{is
convex on}\ (-\infty,0) \nonumber \\
   &\mbox{and concave on}\ (0,\infty)
\end{align}
(this can also be deduced from the formula for $N$ in Theorem 2).

For a class $L$ distribution $\mu=[a,\sigma^2,M]$ on the real line,
Corollary 3 gives that
\begin{align}
\mu \in L^f\  \mbox{iff}\ -x&(dL_{M}(x)/dx)\ \mbox{is convex on}\
(-\infty,0) \nonumber \\
   &\mbox{and concave on}\ (0,\infty).
\end{align}
Some examples of class $L^f$ distributions are provided by the following:
\begin{prop}
If $\eta_1,\eta_2,...$ are i.i.d. Laplace variables (with density
$\frac{1}{2}e^{-|x|})$ and $\sum_{1}^{\infty}a_{k}^2 < \infty$, $a_k >0$,
then $\mu = \mathcal{L}(\sum_{1}^{\infty}a_{k}\eta_k)$ is selfdecomposable
and its background driving distribution $\nu$ is s-selfdecomposable. In
other words, $\mu$ has the factorization property.
\end{prop}
\textit{Proof.} From Jurek (1996) we get that $\mu \in L$ (because Laplace
distributions are selfdecomposable and $L$ is a closed semigroup) and its
L\'evy spectral function has a density of the form: $x \to
\sum_{k=1}^{\infty}exp(-a_k^{-1}|x|)/|x|$. Hence it satisfies the
conditions in (25) and therefore $\mu \in L^f$.

\medskip
1. \textit{L\'evy's stochastic area integrals.} For
$\mathbb{B}_{t}=(B^{1}_t,B^{2}_t)$, Brownian motion on
$\mathbb{R}^{2}$, the process
\[
\mathcal{A}_{t} = \int_0^t B^{1}_sdB^{2}_s-B^{2}_sdB^{1}_s,\ \ t>0,
\]
is called \textit{L\'evy's stochastic area integral}. It is well-known that
for fixed $u>0$, and $a = (\sqrt{u},\sqrt{u})\in \Rset^2 $ we have
\begin{equation}
\chi(t) = E[e^{it\mathcal{A}_{u}}|B_{u}= a] = \frac{tu}{\sinh tu} \cdot
\exp \{-(tu\coth tu-1)\},\quad t\in \Rset
\end{equation}
(cf. L\'evy (1951) or Yor (1992), p.19). Hence the characteristic function
$\chi$ is equal to the product of $\phi(t)  = tu/\sinh tu$, which is
selfdecomposable, and $\psi(t) = \exp[-(tu \coth tu-1)]$, which is its
background driving characteristic function; cf. Jurek (2001), Example B.
However, from Jurek (1996), Example 1, we have that $\phi$ is the
characteristic function of the r.v. $D_1 =
\sum_{k=1}^{\infty}k^{-1}\eta_{k}$, where the $\eta_k$ are as in
Proposition 3. Thus Proposition 3 gives that $\psi$ is the char. f. of an
s-selfdecomposable distribution. Consequently, by Theorem 1, $\chi$ is also
a selfdecomposable char. f.

Using Proposition 3 in Jurek (2001) (or the above relation (26)) we may
find the BDLP for $\chi(t)$.

\medskip
2. \textit{Wenocur integrals.} Let $B_{t}$, $t\in [0,1]$, be
Brownian motion and $Z$ be an independent standard normal random
variable. Then from Wenocur (1986) or Yor (1992), p.19, for a
particular choice of parameters we have that
\begin{equation}
E[e^{itZ\sqrt{\int_{0}^{1}(B_{s}\pm 1)^{2}ds}}]= (\cosh t)^{-1/2}
\cdot \exp(-2^{-1}t\,\tanh t).
\end{equation}
Thus the distribution of $Z\sqrt{\int_{0}^{1}(B_{s}\pm 1)^{2}ds}$
corresponds to the convolution of a class $L$ distribution with
characteristic function $(\cosh t)^{-1/2}$ (which, in fact, is the
characteristic function of a convolution square root of the law of
the r.v. $D_2 = \sum_{k=1}^{\infty}(2k-1)^{-1}\eta_{k}$; cf. Jurek
(1996)) and its BDPD, with char. f. $\exp (-2^{-1}t\tanh t)$,
which is in the class $\mathcal{U}$, by Proposition 3.
Consequently, the above product is also the char. f. of a class
$L$ distribution.

(The fact that in (26) and (27) we have convolutions of selfdecomposable
distributions with their background driving probability distributions was
already observed in Jurek (2001).)

\medskip
3. \textit{Gamma and related distributions.} (a) Let
$\gamma_{\alpha,\lambda}$ be the \textit{gamma distribution} with
probability density $\frac{\lambda^{\alpha}}{\Gamma(\alpha)}x^{\alpha -1}
e^{-\lambda x}1_{(0,\infty)}(x)$. It is easy to see that it is
selfdecomposable (its L\'evy spectral function $L_M$ satisfies the equation
$dL_{M}(x)/dx = \alpha e^{-\lambda x}/x\ \mbox{for}\ x>0$), and its BDRV is
the compound Poisson r.v. $Pois(\alpha \gamma_{1,\lambda}$), i.e., its
jumps are exponentially distributed. Thus by (24), the BDRV is
s-selfdecomposable; in fact it is \textit{s-stable.} The latter are limits
as in (2) where the $\xi_n 's$ are identically distributed; cf. Jurek
(1985), formula (4.3), p. 606.

(b) Let $\rho_{\alpha }$ be the \textit{Bessel distributions} given by the
probability density functions
\begin{equation*}
f_{\alpha }(x)=\exp (-\alpha -x)(x/\alpha )^{(\alpha -1)/2}I_{\alpha -1}(2%
\sqrt{\alpha x}),\text{ }x>0\text{, \ }\alpha >0,
\end{equation*}
where $I_{\alpha -1}(x)$ is \textit{the modified Bessel function}
with index $\alpha -1$. Then $\rho_{\alpha}=\gamma_{\alpha,1}\ast
Pois(\alpha \gamma_{1,1})$. (Iksanov and Jurek (2003) showed that
the Bessel distribution $\rho_{\alpha}$ is a shot-noise
distribution.)

(c) Thorin's distributions (the class $\mathcal{T}$ of generalized gamma
distributions) have the factorization property, as they are obtained from
gammas, their translations, and weak limits; cf. Proposition 1 and case (a)
above, or see Theorem 3.1.1 in Bondesson (1992).

(d) Dufresne (1998) (see p. 295), studied distributional equations
involving symmetrized gamma r.v.'s. These are distributions whose
probability density functions are of a form
\[
p_{\alpha}(x)=\frac{2^{-\alpha +1/2}|x|^{\alpha -1/2}}
{\pi^{1/2}\Gamma(\alpha)}K_{\alpha -1/2}(|x|),\ \ x\in \Rset,
\]
where the $K_{\beta}$ are the MacDonald functions. Since gamma rv's are in
$L$, convergent series of symmetrized gamma rv's provide distributions with
the factorization property.

\medskip
\begin{rem}
Recall that a d.f. $F$ on $\Rset$ is \textit{unimodal with mode at $0$} iff
$F(x) - xF^{\prime}(x)$ is a d.f., or equivalently, iff
$F(x)=\int_{0}^{1}H(x/t)dt$ for some d.f. $H$. Moreover, $H$ may be chosen
to be equal to $F(x)-xF^{\prime }(x)$ a.e. Note that the above relation is
the same as (10) (on the level of L\'evy exponents) or the conditions
described in Theorem 2 (on the level of L\'evy spectral measures).
\end{rem}

\medskip

\noindent Cybernetics Faculty, Kiev T. Shevchenko National
University, 01033 Kiev Ukraine. [E-mail: iksan@unicyb.kiev.ua]

\medskip
\noindent Institute of Mathematics, University of Wroc\l aw,
50-384 Wroc\l aw, Poland. [E-mail: zjjurek@math.uni.wroc.pl]

\medskip
\noindent Department of Mathematics, Wayne State University,
Detroit, MI 48202, USA. [E-mail: berts@math.wayne.edu ]
\end{document}